\documentclass{article}
\usepackage[T1]{fontenc}
\usepackage{lmodern}
\usepackage[latin1]{inputenc}
\usepackage{latexsym}                 
\usepackage{color}                 
\usepackage{amssymb}
\usepackage{amsthm}
\usepackage{graphicx}
\usepackage{hyperref}
\usepackage{verbatim}
\usepackage{mathptmx}      % use Times fonts if available on your TeX system
\def\eref#1{(\ref{#1}%
)}
\def\RSref#1{\ref{#1}%
}
\def\RSlabel#1{\label{#1}%
}
\def\RScite#1{\cite{#1}%
}

\newcommand{\bql}[1]{\refstepcounter{Dummy}%
\begin{equation}\label{#1}%
}

\def\filename#1{}
\newcommand{\eq}{\end{equation}}
\def\fa{\hbox{ for all }}

\def\dfrac#1#2{\displaystyle{\frac{#1}{#2}   }}

\def\b1{\mathbf 1}

\newcommand{\R}{\ensuremath{\mathbb{R}}}
\def\calO{\mathcal{O}}

\def\biglf{\par\bigskip\noindent}
\newtheorem{Dummy}{Dummy}[section]

{\em}
{\em}
\newtheorem{Lemma}[Dummy]{Lemma}{\em}
\newtheorem{Corollary}[Dummy]{Corollary}{\em}
{\em}
\newtheorem{Definition}[Dummy]{Definition}{\em}
\newtheorem{Theorem}[Dummy]{Theorem}{\em}
{\em}
\newtheorem{Example}[Dummy]{Example}{\rm}
{\rm}
{\em}
{\em}
{\em}

\begin{document}
\begin{center}
{\bf All Well--Posed Problems have Uniformly Stable and Convergent Discretizations}
\biglf
Robert Schaback \\
Univ. Göttingen\\
schaback@math.uni-goettingen.de\\
http://num.math.uni-goettingen.de/schaback/research/group.html
\biglf
Draft of \today
\end{center}
{\bf Abstract:}
This paper considers a large class of linear operator equations, including 
linear boundary value problems for partial  
differential equations, and treats them as linear recovery problems for {\em objects}
from
their {\em data}. {\em Well--posedness} of the problem means that this recovery is continuous.
{\em Discretization} recovers restricted {\em trial} objects from restricted
{\em test} data,
and it is well--posed or {\em stable}, if  this {\em restricted} recovery is continuous.
After defining a general framework for these notions, this paper proves that  
all well--posed linear problems have stable and refinable 
computational discretizations with a stability 
that is determined by the well--posedness of the problem
and independent of the computational discretization. The solutions
of discretized problems converge when enlarging the trial spaces, and the convergence rate 
is determined by how well the full data of the object solving the full problem
can be approximated by the full data of the trial objects.
%of the discretization. 
This allows very simple proofs of convergence rates for
generalized finite elements, symmetric and unsymmetric Kansa--type collocation,
and other meshfree methods like  Meshless Local Petrov--Galerkin techniques.  
It is also shown that for a fixed trial space, weak formulations have a slightly
better convergence rate than strong formulations, but at the expense of
numerical integration.
Since convergence rates are reduced to those coming 
from Approximation Theory, and since trial spaces are
arbitrary, this also covers various spectral and pseudospectral methods.
All of this is illustrated by examples.
%%%%%%%%%%%%%%%%%%%%%%%%%%%%%%%%%%%%%%%%%%%%%%%%%%%%%%%%%%%%%%%%%%%%%%%%%%
\section{Introduction and Summary}\RSlabel{Introduction}
This paper focuses on mathematical problems that have solutions $u$ in some 
normed linear space $U$ over $\R$ satisfying infinitely many 
linear conditions that we write as
\bql{eqruf}
\lambda(u)=f_\lambda \fa \lambda\in \Lambda\subseteq U^*
\eq
with given real numbers $f_\lambda$ and continuous linear functionals $\lambda$
on $U$ collected into a set $\Lambda\subset U^*$. We call the real numbers
$\{f_\lambda\}_{\lambda\in \Lambda}$ the {\em data} that hopefully allow to identify the
{\em object} $u$, which will in many cases be a multivariate function on some
domain.
Solving \eref{eqruf} for $u$ from given data $\{f_\lambda\}_{\lambda\in
  \Lambda}$ 
is a {\em recovery problem}, and we view it as posed in an abstract mathematical
setting
that is not directly accessible for computation. It can be called an {\em
  analytical}
problem in contrast to the {\em computational} problems that will follow later.
The transition from an analytical problem to a computational problem will be
called {\em discretization}.   
\biglf
Typical special cases arise when solving partial differential equations (PDEs).
The object $u$ to be recovered is always an element of some space $U$   
of real--valued functions on a domain $\Omega$, but weak and strong formulations 
of PDEs use very different types of {\em data}, namely either integrals against
test functions or derivative values at evaluation points, plus boundary
conditions of various forms. If two problems use
different data to identify the same object, we consider them as different here.
\biglf
The PDE applications of \eref{eqruf} suggest to view the application
of the functionals $\lambda\in \Lambda$ as {\em testing} a {\em trial} object
$u$. {\em Discretization} will then fix a finite--dimensional {\em trial subspace}
$U_M\subset U$ and a finite set $\Lambda_N$ of {\em test functionals} from $\Lambda$. 
We pursue this distinction between the trial and the test side of \eref{eqruf}
throughout this paper.
%, starting in Section \RSref{SecIntroD} with
%definitions of {\em data}, {\em discretization}, and {\em refinement}.   
\biglf
If a problem in {\em Mathematical} Analysis is well--posed, it should have a
discretization in {\em Numerical} Analysis that is also well-posed. This requires
to derive some sort of numerical stability of well--chosen discretizations
from the well--posedness property of the underlying analytical problem.
This paper proves the above statement under mild additional assumptions
after stating clearly in Sections \RSref{SecIntroD} 
and \RSref{SecIntroDD} what is to be understood
by {\em well--posedness} of a {\em problem} and its {\em discretization}.  
It turns out in Section  \RSref{SecIntroWPDP} that one can choose refinable discretizations 
that have stability properties depending only on the well--posedness
of the given problem, not on the discretizations chosen. This depends 
crucially on what we call a {\em monotone refinable dense (MRD)} discretization 
strategy in Section \RSref{SecIntroDD}. 
In an older discretization theory
\RScite{schaback:2007-1,schaback:2010-2,boehmer-schaback:2013-1},
error bounds and convergence results depended on {\em stability inequalities}
that needed complicated proofs \RScite{rieger-et-al:2010-1},
while this paper shows that one can always
enforce uniform stability  by sufficiently thorough testing. 
\biglf
The resulting discretized linear problems will be overdetermined
due to this stabilization, and should be solved approximately by 
minimizing residuals. Section \RSref{SecIntroSDP} deals with  
this, and shows that the final error bounds and convergence rates
are determined by how well the data of the true solution
can be approximated by the data of elements of the trial space.
We call this {\em Trial Space Data Approximation}. 
In particular, error bounds and convergence rates are independent of the 
details of testing.
\biglf
Section \RSref{SecIntroND} extends the previous results to ill--posed problems
and noisy data, while Sections \RSref{SecDiHS} and \RSref{SecORiHS}
specialize to recovery in Hilbert spaces, where uniformly stable 
and sometimes optimal discretizations are readily available. These generalize 
finite elements and symmetric collocation, as will be explained
in Section \RSref{SecSDP} when it comes to examples.
\biglf
Stability can be spoiled by bad bases. Therefore this paper ignores bases
and focuses on spaces instead, up to Section \RSref{SecB} where the
influence of bases on the trial and the test side is studied. 
A very common class of bases are the {\em nodal} bases used in classical
piecewise linear finite elements and  various meshless methods.
Many application papers report good stability properties of these,
and Section \RSref{SecNod} provides a 
fairly general mathematical proof, showing that convergence in the nodes can be derived from
convergence of the {\em Trial Space Data Approximation}.
\biglf
The paper closes with a number of examples that apply the above theory. 
Polynomial interpolation in Section \RSref{SubSecInt} illustrates  that the
stabilization results of this paper imply quite some {\em overtesting} , i.e. oversampling 
on the test side to guarantee uniform stability. Furthermore, it 
points out how
spectral methods are covered and why weak formulations yield slightly
faster convergence than strong formulations, 
though for weaker norms and at the expense of numerical integration.
\biglf 
Section \RSref{SecSDP} deals with the standard setting for 
finite elements for homogeneous boundary conditions, 
showing that it fits perfectly into the framework,
including extension to other trial spaces and a Petrov--Galerkin treatment.
\biglf
The remaining examples address the standard Poisson problem with Dirichlet
boundary conditions, for simplicity.
Section \RSref{SubSecCM} focuses
on collocation as a typical strong problem. This covers
various kinds of meshless methods, including
Kansa's   unsymmetric collocation, and it is shown how to derive 
specific convergence rates depending on the trial spaces chosen.
The weak Dirichlet case is handled  in Section
\RSref{SecWDP}, and a comparison of convergence rates for the strong and weak
formulations, using the same trial spaces, is provided in Section \RSref{SecWSC}. 
\biglf
Finally, Atluri's Meshless Local Petrov Galerkin (MLPG) scheme
\RScite{atluri:2005-1}
is treated in
Section \RSref{SecMLPG}. This includes error bounds and convergence rates
for different variations of the method,
but it  was necessary to include a first proof of well--posedness 
of the local weak form behind MLPG.
\biglf
Summarizing, this paper shows that under mild hypotheses 
\begin{enumerate}
 \item all well--posed problems have uniformly stable discretizations
made possible by sufficiently extensive testing, and
\item convergence rates for such discretizations can be played back to
known convergence rates of {\em Trial Space Data Approximation},
i.e. the approximation of the data of the true solution   by
the data of the trial elements. These rates depend on what ``data'' means 
and are taken in the norm arising in the well--posedness condition.
\item Weak and strong formulations of a given background problem will 
have different definitions of ``data'' and will need different 
versions of well--posedness, and these differences enter into the previous item 
and influence the convergence rates, even when trial spaces are the same for both
formulations.
\item For a given fixed trial space, it is shown that in standard applications the
  weak 
formulation converges slightly faster than  the strong formulation. 
\item Nodal bases have a stability advantage over other bases.
\end{enumerate}   
On the downside, the test strategies guaranteeing 
uniform stability are only shown to exist, they are not constructed.
Future work needs  explicit sufficient conditions on 
specific test strategies to guarantee uniform stability.
This can be done by {\em greedy testing} as touched 
in earlier papers on adaptivity 
\RScite{schaback-wendland:2000-2,hon-et-al:2003-1,ling-schaback:2004-1,schaback:2013-2}.  
Finally, emphasis so far is only on errors, convergence rates, and stability of
algorithms,
but not on computational efficiency. It is a major challenge to 
relate the achievable convergence rates and stability properties to
computational efficiency. Anyway, this paper provides a starting point
towards these goals.
%%%%%%%%%%%%%%%%%%%%%%%%%%%%%%%%%%%%%%%%%%%%%%%%%%%%%%%%%%%%%%%%%%%%%%%%%%
\section{Problems, Data, and Well--Posedness}\RSlabel{SecIntroD}
Behind \eref{eqruf} there is a linear {\em data map} $D\;:\;U\mapsto \R^\Lambda=:V$ 
that takes each $u\in U$
into the set of values $\{\lambda(u)\}_{\lambda\in \Lambda}$ in the {\em data
  space} $V$. This allows to rewrite \eref{eqruf} as an operator equation
\bql{eqDuf}
D(u)=f
\eq
for some given $f$ in the data space $V$. Each operator equation can be
formally interpreted this way, e.g. by defining $\Lambda$ as the set
of all functionals $\mu\circ D$ when $\mu$ varies in the unit sphere of $V^*$.
\begin{Example}\RSlabel{ExaPP}{\rm
As an illustration, consider the standard Dirichlet problem 
$$
\begin{array}{rcll}
-\Delta u &=& f & \hbox{ in } \Omega\subset\R^d\\ 
u &=& g & \hbox{ in } \Gamma:=\partial \Omega 
\end{array}
$$
where $f$ and $g$ are given functions on $\Omega$ and $\Gamma$.
A {\em strong} formulation writes it in the form \eref{eqruf} with functionals
$$
\begin{array}{rcll}
\lambda_x &:& u\mapsto -\Delta u(x),& x\in \Omega\\
\mu_y &:& u\mapsto u(y),& y\in \Gamma\\
\end{array}
$$
on some space $U$ where both types of functionals are continuous.
The standard FEM algorithms use {\em weak} functionals
$$
\lambda_v\;:\;u\mapsto (\nabla u,\nabla v)_{L_2(\Omega)}\,\fa v\in H_0^1(\Omega)
$$ 
and add the functionals $\mu_y$ for points on the boundary. We postpone further
details to sections \RSref{SubSecCM} and \RSref{SecWDP}, but remark that 
the data maps differ considerably.}\qed
\end{Example}  
We give the data space a norm structure by requiring that 
\bql{eqDuVlamu} 
\|Du\|_V:=\sup_{\lambda\in \Lambda}|\lambda(u)| \fa u\in U
\eq
is a norm on $D(U)$ that we assume to be extended to $V$, if not $V=D(U)$. 
We shall call this the {\em data norm}, and note that it leads to a seminorm
\bql{equDDuV}
\|u\|_D:=\|Du\|_V=\sup_{\lambda\in \Lambda}|\lambda(u)| \fa u\in U
\eq
on the object space $U$.
This is well--defined if 
all functionals in $\Lambda$ are uniformly bounded. 
We assume existence of the data norm from now on, but remind the reader that
renormalization of functionals changes the data norm and all issues 
depending on it, like the well--posedness conditions that we introduce later.
\begin{Definition}\RSlabel{DefAnPro}
An {\em analytic problem} in the sense of this paper consists of
\begin{enumerate}
 \item a linear normed {\em object} space $U$,
\item a set $\Lambda$ of linear functionals on $U$ leading to a {\em data map} $D$ as
  in \eref{eqDuf} 
\item with values in a normed {\em data space} $V$  such that 
\item \eref{eqDuVlamu} holds and is a norm on $V$.
\end{enumerate} 
\end{Definition} 
Unique solvability of the problem \eref{eqruf} or \eref{eqDuf} 
requires that $u\in U$ vanishes if all data 
$\lambda(u)$ for all $\lambda\in \Lambda$ vanish, or that  $D$ is injective,
or that $\|.\|_D$ is a norm.
A somewhat stronger 
and quantitative notion is {\em well--posedness}:
\begin{Definition}\RSlabel{DefWellPos}
An analytic problem in the sense of Definition \RSref{DefAnPro}
is {\em well--posed with respect to a  well--posedness norm} $\|.\|_{WP}$ on $U$
if there is a constant $C$ such that 
a {\em well--posedness
  inequality} 
\bql{equWPCDu}
\|u\|_{WP}\leq C\|Du\|_V=C\|u\|_D \fa u\in U
\eq
holds. 
\end{Definition} 
This  means that
$D^{-1}$ is continuous as a map $D(U)\to U$ in the norms $\|.\|_V$ and
$\|.\|_{WP}$. The well--posedness norm $\|.\|_{WP}$ on $U$ will often be weaker
than the norm $\|.\|_U$ on $U$ needed to let the data be well--defined.
\biglf
In the context of Example \RSref{ExaPP}, the strong problem
leads to well--posedness with $\|.\|_{WP}=\|.\|_{\infty,\Omega}$, while the weak form has
$\|.\|_{WP}=\|.\|_{L_2(\Omega)}$. Details will follow in \RSref{SubSecCM} and
\RSref{SecWDP}, but we remark here that deriving 
computationally useful well--posedness inequalities is a
serious issue that is not satisfactorily addressed by theoreticians,
because they do not use computationally useful norms on the data space.  
For instance, the continuous dependence of solutions of elliptic problems
on the boundary data is often expressed by taking Sobolev trace spaces
of fractional order on the boundary, and these spaces are far from being accessible
for computation. The examples will shed some light on this issue.
\biglf
Future research in Applied Mathematics should target 
practically useful well--posedness results based on norms that are closer to
computation.    
%%%%%%%%%%%%%%%%%%%%%%%%%%%%%%%%%%%%%%%%%%%%%%%%%%%%%%%%%%%%%%%%%%%%%%%%%%
\section{Trial Space Data Approximation}\RSlabel{SecIntroTSA}
We now perform the first step of {\em discretization} by choosing a 
finite--dimensional {\em trial
  space} $U_M\subset U$. This allows us to approximate the data $D(u^*)\in V$
by data $D(u_M)$ for all trial elements $u_M\in U_M$ in the data norm $\|.\|_V$, 
and we denote the best approximation by $u_M^*$, i.e.
\bql{eqDuDuMV}
\|Du^*-Du_M^*\|_V=\displaystyle{\min_{u_M\in U_M}\|Du^*-Du_M\|_V.   } 
\eq
We shall rely on Approximation Theory to provide upper bounds for this, and
for convergence rates for $\|Du^*-Du_M^*\|_V\to 0$ 
if the spaces $U_M$ get larger and larger. These rates will
crucially
depend on the smoothness of $u^*$, the trial spaces $U_M$, and the data map $D$. 
For trial spaces in spectral methods, these convergence rates may be very large,
and there may even be exponential convergence. We call \eref{eqDuDuMV} the
{\em Trial Space Data Approximation}, but we keep in mind that the above 
approximation problem is computationally hazardous, because it involves
infinitely many data. We can assess $u_M^*$ only in theory, not in practice.
\biglf
If the problem is well--posed in the sense of Definition
\RSref{DefWellPos}, the 
error bounds and convergence rates 
of the {\em Trial Space Data Approximation} immediately carry over to 
error bounds and convergence rates in the well--posedness norm, via
$$
\|u^*-u_M^*\|_{WP}\leq C \|Du^*-Du_M^*\|_V,
$$ 
and independent of the chosen trial space. This means that Approximation
Theory provides convergence rates for certain approximate solutions of 
certain well--posed analytic problems, but these approximate solutions are
computationally inaccessible. 
\biglf
In the context of Example \RSref{ExaPP}, the functions of the trial space 
have to approximate function values on the boundary in both the strong and the
weak case. But for the strong form we have to approximate second derivatives,
while
the weak form only has to approximate first derivatives. Furthermore,
the well--posedness norms are different. This will lead to
different convergence rates in Section \RSref{SecWSC}.  
%%%%%%%%%%%%%%%%%%%%%%%%%%%%%%%%%%%%%%%%%%%%%%%%%%%%%%%%%%%%%%%%%%%%%%%%%%
\section{MRD Discretizations}\RSlabel{SecIntroDD}
\biglf
In what follows, we shall show how to 
discretize the test side of an analytic problem in the sense of Definition
\RSref{DefAnPro} in such a way
that a uniformly stable and finite computational strategy exists that
provides approximations $\tilde u_M\in U_M$ with   
$$
\|u^*-\tilde u_M\|_{WP}\leq 2\|u^*-u_M^*\|_{WP}\leq 2C \|Du^*-Du_M^*\|_V.
$$ 
This implies that Approximation
Theory provides convergence rates for certain {\em finitely and stably computable }
approximate solutions of 
certain well--posed analytic problems. The convergence will take place in $U$
under the
well--posedness norm $\|.\|_{WP}$, and the convergence rate will be the
convergence rate of the {\em Trial Space Data Approximation}.
Our main tool will be a 
{\em monotonic refinable dense (MRD)} discretization of the data space $V$
that we describe now. 
\biglf
No matter what the data {\em map} is, the data {\em space}  $V$ should allow 
some form of {\em discretization} for computational purposes. We model this by
{\em restriction} maps 
$$
R_N\;:\;\{f_\lambda\}_{\lambda\in \Lambda} \mapsto \{f_\lambda\}_{\lambda\in
  \Lambda_N}\in V_N=\R^{|\Lambda_N|}
$$ 
that map $V$ into finite--dimensional {\em data spaces}
$V_N$ over $\R$. The discretizations use   
{\em restricted data} belonging to finite subsets   $\Lambda_N$ of $\Lambda$,
and these  data enter practical computation. 
\biglf
On the spaces $V_N$ we introduce the norm  
$$
\|R_N\{f_\lambda\}_{\lambda\in \Lambda}\|_{V_N}=\|\{f_\lambda\}_{\lambda\in
  \Lambda_N}\|_{V_N}
:=\displaystyle{\max_{\lambda\in \Lambda_N}|f_\lambda|   }  
$$
and we get the  {\em monotonicity} property 
$$
 \|R_Mv\|_{V_M} \leq \|R_Nv\|_{V_N} \fa v\in V \hbox{ and all } \Lambda_M\subseteq \Lambda_N.
$$ 
{\em Refinement} of two discretizations defined by sets $\Lambda_M$ and
$\Lambda_N$  works by taking $\Lambda_M\cup \Lambda_N$, and by the monotonicity
property this
will weakly increase the discrete norms. Finally, we have  
\bql{eqvVDens}
\|v\|_V:=\displaystyle{\sup_{R_N, V_N} \|R_Nv\|_{V_N}   }\fa v\in V, 
\eq
following from \eref{eqDuVlamu}.  
\biglf
But there are applications where restrictions are not defined by taking 
{\em all
possible} finite 
subsets of functionals. They might require background triangulations, e.g. for
finite elements, and their refinement does not simply involve taking a union
of two finite sets of functionals. We can generalize the above notions
by ignoring functionals:
\begin{Definition}\RSlabel{DefMRD}
An {\em MRD discretization} of a data space $V$ consist of a set of
{\em restrictions} $(R_N,V_N)$ with the properties  
\begin{enumerate}
\item $V_N$ is a normed linear space with $\dim V_N<\infty$ and norm $\|.\|_{V_N}$,
\item $R_N\;:V\to V_N$ is linear,
\item there is a partially defined {\em refinement relation} $\preceq$ on the
  restrictions such that
\item  $(R_M,V_M) \preceq (R_N,V_N)$ implies  $\|R_Mv\|_{V_M} \leq \|R_Nv\|_{V_N} \fa v\in V$,
\item  for each two admissible restrictions $(R_M,V_M),\;(R_N,V_N)$ there is a restriction
$(R_P,V_P)$ such that $(R_M,V_M) \preceq (R_P,V_P)$  and $(R_M,V_M) \preceq
  (R_P,V_P)$,
\item \eref{eqvVDens} is a norm, when the sup is taken over all admissible restrictions.
\end{enumerate} 
\end{Definition} 
This axiomatic framework is open for further discussion, of course, but we assume it
in what follows. We refer to the last three properties as {\em monotonicity},
{\em refinement}, and {\em density}, using the term {\em MRD discretization} 
for all six properties. Note that the norm arising in the density property must
be the data norm that is used in the well--posedness inequality
\eref{equWPCDu}. 
\biglf
The discussion preceding Definition \RSref{DefMRD} proved
\begin{Theorem}\RSlabel{TheExMRD}
Each analytical problem of the form  \eref{eqruf} in the sense of Definition \RSref{DefAnPro}
  %with a data norm \eref{eqDuVlamu} 
has a MRD discretization via taking finite subsets of functionals. \qed
\end{Theorem} 
For Example \RSref{ExaPP}, it is clear that one can focus on finitely many
functionals when it comes to finite computations, but it is by no means clear 
which and how many are to be taken to allow a uniformly stable computational
method.
The refinement in the FEM case is not quite standard, but will satisfy
Definition \RSref{DefMRD}, because it still uses finite subsets of functionals.   
%%%%%%%%%%%%%%%%%%%%%%%%%%%%%%%%%%%%%%%%%%%%%%%%%%%%%%%%%%%%%%%%%%%%%%%%%%
\section{Well--Posedness of Discretized Problems}\RSlabel{SecIntroWPDP}
If we use a MRD restriction $(R_N,V_N)$ on the data together with a chosen trial
space $U_M$, we can 
pose the {\em discretized problem } as the linear system
\bql{eqRNDUM}
R_NDu_M=R_NDu^*
\eq
to be solved for $u_M\in U_M$, where the computational input 
data are provided by the restriction $R_NDu^*$ of the data of an exact solution
$u^*$. Such systems will usually be overdetermined. 
\biglf
Since the well--posedness condition 
\eref{equWPCDu} also holds on the trial space, the discretized problem \eref{eqRNDUM}
is automatically well--posed or {\em stable} in the sense
$$
\|u_M\|_{WP}\leq C(U_M, V_N) \|R_NDu_M\|_{V_N} \fa u_M\in U_M
$$
if  we can prove
\bql{equMDDuMV}
%\|u_M\|_D=
\|Du_M\|_V\leq C(U_M, V_N) \|R_NDu_M\|_{V_N} \fa u_M\in U_M
\eq
for some {\em stability constant} $C(U_M,V_N)$. 
\biglf
We now can state our central result, to be
proven later in a somewhat more general form.
\begin{Theorem}\RSlabel{TheGenStab}
Assume an analytic  problem \eref{eqDuf} with an MRD discretization.
If $U_M$ is an arbitrary finite--dimensional subspace of $U$, there always is
a restriction $(R_N,V_N)$ such that
\bql{equMD2RNDuM}
\|u_M\|_{D}\leq 2 \|R_ND(u_M)\|_{V_N}\fa u_M \in U_M.
\eq
This holds without assuming well--posedness. If the latter is assumed by \eref{equWPCDu}, we have 
\bql{equMU2C}
\|u_M\|_{WP}\leq 2C \|R_ND(u_M)\|_{V_N}\fa u_M \in U_M
\eq
with the constant $C$ from \eref{equWPCDu}.
\end{Theorem} 
In contrast to \eref{equMDDuMV}, the constants in \eref{equMD2RNDuM} and
\eref{equMU2C}
are independent of $U_M$ and $V_M$, proving a {\em uniform} well--posedness
or stability of the discretized problem for a
rather sensible choice of $V_N$ after an arbitrary selection of $U_M$.
Section \RSref{SubSecInt} will show that this uniformity may require some hidden amount of
{\em oversampling},
i.e. the dimension of $V_N$ may be much larger than the dimension of $U_M$.
We call this  {\em overtesting}, because it occurs on the test
side of the problem. Theorem \RSref{TheGenStab} does not give any practical hints how to 
care for uniformly stable testing, it just proves existence. The necessary 
amount of overtesting to achieve uniformly stability is left open. 
\biglf
It is a common observation that many instabilities arise from badly chosen
bases. They sometimes disappear after introduction of better bases. To 
identify instabilities that can be blamed to bad bases, 
we refrain from introducing bases as far as possible in this paper, focusing on
spaces instead of bases.
%\biglf
%We can rephrase the above theorem by stating that all well--posed analytic
%problems have uniformly stable discretizations, provided that the data space $V$
%has an MRD discretization. The trial spaces can be chosen
%arbitrarily, but the testing must then be done with sufficiently many well--chosen
%conditions to achieve {\em uniform} well--posedness of the discretized
%problems. 
%It handles the situation where a trial space $U_M$
%is chosen first, followed by a good choice of a restriction $R_N$. On can also
%start with a restriction $R_N$ and look for good trial spaces $U_M$ afterwards.
%We shall pursue this approach later \red{Wo????} as well. 
%
%\red{Greedy testing?} 
%%%%%%%%%%%%%%%%%%%%%%%%%%%%%%%%%%%%%%%%%%%%%%%%%%%%%%%%%%%%%%%%%%%%%%%%%%
\section{Well--Posedness of Data Discretizations}\RSlabel{SecIntroWPDD}
Inspection of \eref{equMD2RNDuM} shows that the analytic problem and its
well--posedness are not relevant for \eref{equMD2RNDuM}, because the
actual well--posedness condition \eref{equWPCDu} enters only into the trivial
transition from \eref{equMD2RNDuM} to \eref{equMU2C}. In fact, everything
follows
already from the notion of a MRD discretization. 
Well--posedness is a later add--on.
\begin{Lemma}\RSlabel{LemDatStab}
Consider a data space $V$ and associated MRD restrictions $(R_N,V_N)$ satisfying the
assumptions of Section \RSref{SecIntroDD}. Then for each finite--dimensional
subspace $W_M$ of $V$ there always is
a restriction $(R_N,V_N)$ such that
$$
\|w_M\|_{V}\leq 2 \|R_Nw_M\|_{V_N}\fa w_M \in W_M.
$$   
\end{Lemma} 
{\bf Proof}:  Define
$K\subset W_M$ as the unit sphere of $W_M$ defined via the norm $\|.\|_V$. 
By compactness,  
for each $\epsilon>0$ we can cover $K$ by finitely many
$\epsilon$--neighborhoods
$$
U_\epsilon(y_j):=\{y\in K\;:\;\|y-y_j\|_V\leq \epsilon  \},\,1\leq j\leq n
$$
with elements $y_1,\ldots,y_n\in K$.  By the density property
\eref{eqvVDens} we can find restrictions
$R_{N_1},\ldots,R_{N_n}$ with associated spaces
$V_{N_1},\ldots,V_{N_n}$ such that  
$$
\|y_j\|_V\leq \|R_{N_j}y_j\|_{V_{N_j}}+\epsilon,\,1\leq j\leq n
$$
and by repeated application of the refinement property 
we can define $R_N$ and $V_N$ as the ``union'' of these, and then 
$$
\|R_{N_j}v\|_{V_{N_j}}\leq \|R_{N}v\|_{V_N}\fa v\in V,\,1\leq j\leq n
$$
by monotonicity.
\biglf
We now take an arbitrary $w_M\in K$ and get some $j,\,1\leq
j\leq n$ with
$\|w_M-y_j\|_V\leq \epsilon$ via the covering. This implies
$\|R_Nw_M-R_Ny_j\|_{V_N}\leq \epsilon$ by the density property, and then
$$
\begin{array}{rcl}
\|R_Nw_M\|_{V_N}
&\geq &
\|R_Ny_j\|_{V_N}-\epsilon\\
&\geq &
\|R_{N_j}y_j\|_{V_{N_j}}-\epsilon\\
&\geq &
\|y_j\|_{V}-2\epsilon\\
&\geq &
\|w_M\|_{V}-3\epsilon\\
&=&
1-3\epsilon\\
\end{array}
$$
proving
$$
\|R_Nw_M\|_{V_N}\geq (1-3\epsilon)\|w_M\|_V
$$
for all $w_M\in W_M$, and the assertion follows for $\epsilon=1/6$.\qed
\biglf
The proof of Theorem \RSref{TheGenStab} now follows by setting $W_M=D(U_M)$ 
with an arbitrary data map $D$.\qed
%%%%%%%%%%%%%%%%%%%%%%%%%%%%%%%%%%%%%%%%%%%%%%%%%%%%%%%%%%%%%%%%%%%%%%%%%%
\section{Solving Discretized Problems}\RSlabel{SecIntroSDP}
After choosing a trial space $U_M$ and getting a suitable data restriction $(R_N,V_N)$
for Theorem \RSref{TheGenStab}, the discretized recovery problem \eref{eqRNDUM}
requires computation of some $u_M\in U_M$ from the data $R_NDu^*$, where $u^*$ 
is the true solution to the analytical problem. This will usually lead to an
overdetermined linear system after choosing bases, but we do not want to
consider bases unless absolutely necessary.
\biglf
The simplest basis--free computational method we could propose is to minimize the residual norm
$\|R_ND(u^*-u_M)\|_{V_N}$ 
over all $u_M\in U_M$, which is a finite--dimensional approximation
problem. 
A good candidate in the trial
space $U_M$ is the best approximation $u_M^*$ to the solution $u^*$, and we are
free to consider the best approximation in one of the norms
$\|.\|_U,\;\|.\|_{WP}$ or $\|.\|_D$. We do not insist here on using the data
norm
and a {\em best} approximation, as
we did in \eref{eqDuDuMV}. Instead, we keep  
the choice of $u_M^*$ free
and 
are satisfied with computing an element $\tilde u_M\in U_M$ with 
\bql{equuM}
\|R_ND(u^*-\tilde u_M)\|_{V_N} \leq C_A\|R_ND(u^*-u_M^*)\|_{V_N},
\eq
however it is calculated, with a fixed constant $C_A\geq 1$ that makes
computational life easier when chosen not too close to one.  We call   $u_M^*$ 
 a {\em comparison trial object}. It 
is usually provided by some result of
Approximation Theory that yields a 
useful bound on the right--hand side of \eref{equuM}. Due to
the monotonicity and density properties, we also have 
$$
\|R_ND(u^*-\tilde u_M)\|_{V_N} \leq C_A\|R_ND(u^*-u_M^*)\|_{V_N}\leq C_A\|D(u^*-u_M^*)\|_{V},
$$
such that the best {\em Trial Space Data Approximation} is always an upper bound.
\biglf
Anyway, \eref{equuM} implies
$$
\begin{array}{rcl}
\|\tilde u_M-u_M^*\|_{D}
&\leq &
2\|R_ND(\tilde u_M-u_M^*)\|_{V_N}\\
&\leq &
2\|R_ND(\tilde u_M-u^*)\|_{V_N}+2\|R_ND(u^*-u_M^*)\|_{V_N}\\
&\leq &
(2C_A+2)\|R_ND(u^*-u_M^*)\|_{V_N}\\
&\leq &
(2C_A+2)\|D(u^*-u_M^*)\|_{V}\\
&=&
(2C_A+2)\|u^*-u_M^*\|_{D}\\
\end{array}
$$
and
$$
\begin{array}{rcl}
\|\tilde u_M-u^*\|_{D}
&\leq&
\|\tilde u_M-u_M^*\|_{D}+\|u_M^*-u^*\|_{D}\\
&\leq&
(2C_A+3)\|u_M^*-u^*\|_{D},
\end{array} 
$$
proving that the error of the computational solution $\tilde u_M$
is up to a factor the same as the error of the comparison trial object $u_M^*$, 
evaluated in the data norm. 
\begin{Theorem}\RSlabel{TheGenConv}
Assume an MRD % uniformly well--posed 
discretization of an % a well--posed 
analytic problem along the lines of the
previous sections.
Then each computational technique to solve the discretized problem
approximatively by an element $\tilde u_M\in U_M$ such that \eref{equuM} holds,
will also guarantee
$$
\|\tilde u_M-u^*\|_{D}\leq (2C_A+3) \|\|u_M^*-u^*\|_{D} 
$$
for any  comparison trial object $u_M^*$.
\end{Theorem} 
\begin{Corollary}\RSlabel{CorConv}
Adding well--posedness to Theorem \RSref{TheGenConv} yields 
$$
\|u^*-\tilde u_M\|_{WP}\leq C(2C_A+3)\|u_M^*-u^*\|_{D}
$$ 
proving that convergence rates in the data norm 
transfer to the same convergence rates in the norm $\|.\|_{WP}$  on $U$
arising in the well--posedness condition 
\eref{eqDuVlamu}.\qed
\end{Corollary}
We summarize what we have so far, for easy reference in the examples.
\begin{Theorem}\RSlabel{TheConv2}
Assume a well-posed analytic problem with an MRD discretization
as in Definitions \RSref{DefAnPro}, \RSref{DefWellPos}, and \RSref{DefMRD}.
Then for arbitrary trial spaces one can choose
uniformly stable test discretizations to get uniformly stable computational methods
based on some form of residual minimization. The convergence rates,
measured in the well--posedness norm, are given by the 
convergence rate of the {\em Trial Space Data Approximation}, i.e. the  
rate in which the
data of the true solution are  approximated by 
the data of comparison trial objects, measured in the data norm. \qed
\end{Theorem} 
This will be applied in the following way. First, one assumes
additional 
regularity of the solution object $u^*$ and fixes 
a well--known approximation process in $U$ that provides good 
comparison trial objects $u_M^*$ for these trial spaces, and
with a very good convergence rate $u_M^*\to u^*$ that may even be spectral 
in a weak norm like $\|.\|_{WP}$. Then these approximations are used for
comparison in the above theory, and the convergence rate in the data norm is
calculated from what is known about the approximation process. Then
we know that this rate is the one that arises when solving the analytic problem,
and it arises in the well--posedness norm.
This may even yield spectral convergence, and we shall provide examples.
But note that the rate of convergence of our discretized solutions of the analytical
problems  is only the rate the convergence obtained {\em after} the data map is
applied, and it involves the norm $\|.\|_{WP}$ used in the well--posedness
condition.
\biglf
The above approach applies to a large variety of well--posed analytic problems,
and shows that for properly chosen  scales of trial spaces $U_M$ and
properly chosen test strategies depending on each $U_M$ one gets 
uniformly stable and convergent computational 
methods with convergence rates that can be derived from results of Approximation
Theory. These rates normally improve with the smoothness of the true solution,
but they also depend on the
data map and the well--posedness norm. For a given PDE problem like in Example
\RSref{ExaPP}, the convergence rates of strong and weak formulations  will be
different, even if the trial spaces are the same. This is due to the fact that
the data maps, data norms, and
well--posedness norms are different. Details will follow in Section
\RSref{SecWSC}.
\biglf
If the true solution necessarily has certain singularities
of a known type, like in elliptic PDE problems
on domains with incoming corners, one should always add 
the correct singular functions to
the trial space. Then the approximation quality of the singular solution in the
augmented
trial space is the same as  the approximation quality of a regular solution in the
original trial space, and this quality will improve with the smoothness
of the regular solution. In this sense, going over to extended trial spaces 
like in the XFEM or GFEM does not need a new theory here.
%%%%%%%%%%%%%%%%%%%%%%%%%%%%%%%%%%%%%%%%%%%%%%%%%%%%%%%%%%%%%%%%%%%%%%%%%%
\section{Noisy Data}\RSlabel{SecIntroND}
Corollary \RSref{CorConv} showed that ill--posed problems can be treated,
if one is satisfied with reproducing the data well. But so far we always have
assumed that the input data are given exactly as data of an existing solution.
If this is not true, a few changes are necessary. We assume that the data map
$D$ is always exact, but the input data for computations are assumed to be polluted by either
noise or errors in evaluating the data functionals. 
This also handles the error committed
by numerical integration when the data functionals of weak PDE problems are 
considered. 
\biglf
The data now consist of a general element $v^*$ of the data space $V$, and we assume
that 
there is an object $u^*\in U$ such that $\|D(u^*)-v^*\|_V$ is small, and we want
to recover this object well, or others with a similarly good data reproduction.
We choose a trial space $U_M\subset U$ as before, and we define
$W_M:=D(U_M)\subset V$
as at the end of Section \RSref{SecIntroWPDD}.  
Clearly, \eref{equuM} now has to be replaced by
\bql{equuMv}
\|R_N v^*-R_ND\tilde u_M\|_{V_N} \leq C_A\|R_N v^*-R_NDu_M^*\|_{V_N},
\eq
because there are no other data at hand. Then  
\begin{Theorem}\RSlabel{TheGenConvNoise}
Assume an MRD % uniformly well--posed 
discretization of an % a well--posed 
analytic problem along the lines of the
previous sections.
Then each computational technique to solve the discretized problem
approximatively by an element $\tilde u_M\in U_M$ such that \eref{equuMv} holds,
will guarantee
\bql{equMu2C3}
\|\tilde u_M-u^*\|_{D}\leq (2C_A+3) \|\|u_M^*-u^*\|_{D}+ (2C_A+2)\|v^*-Du^*\|_{V}
\eq
for any $u^*\in U$.
\end{Theorem}
{\bf Proof}: 
We proceed like above, via
$$
\begin{array}{rcl}
\|\tilde u_M-u_M^*\|_{D}
&\leq &
2\|R_ND(\tilde u_M-u_M^*)\|_{V_N}\\
&\leq &
2\|R_ND\tilde u_M-R_Nv^*\|_{V_N}+2\|R_Nv^*-R_NDu_M^*\|_{V_N}\\
&\leq &
(2C_A+2)\|R_Nv^*-R_NDu_M^*)\|_{V_N}\\
&\leq &
(2C_A+2)\left( \|R_Nv^*-R_NDu^*\|_{V_N}+ \|R_NDu^*-R_NDu_M^*\|_{V_N} \right)\\
&\leq &
(2C_A+2)\left( \|v^*-Du^*)\|_{V}+ \|D(u^*-u_M^*)\|_{V}\right)\\
&\leq &
(2C_A+2)\left( \|v^*-Du^*)\|_{V}+\|u^*-u_M^*\|_{D}\right)\\
\end{array}
$$ 
and get
$$
\begin{array}{rcl}
\|\tilde u_M-u^*\|_{D}
&\leq&
\|\tilde u_M-u_M^*\|_{D}+\|u_M^*-u^*\|_{D}\\
&\leq&
(2C_A+3)\|u_M^*-u^*\|_{D}+ (2C_A+2)\|v^*-Du^*\|_{V}.\qed
\end{array} 
$$ 
The inequality \eref{equMu2C3} shows that errors in the data functionals, e.g.
integration errors for weak data, can spoil the convergence unless they
are at least as small as the error committed by the comparison object $u_M^*$
in the data norm. For trial spaces that allow  fast convergence, the 
admissible errors in the data functionals are severely restricted by this observation.
\biglf
One can go into \eref{equuMv} by choosing $u_M^*$ as the minimizer of
$\|v^*-D(u_M)\|_V$ over all $u_M\in U_M$. Then
$$
\|R_Nv^*-R_ND(\tilde u_M)\|_{V_N} \leq C_A\|R_Nv^*-R_ND(u_M^*)\|_{V_N}
\leq C_A\|v^*-Du_M^*\|_{V},
$$
and the proof of Theorem \RSref{TheGenConvNoise} yields
$$
\begin{array}{rcl}
\|\tilde u_M-u_M^*\|_{D}
&\leq &
(2C_A+2)\|R_Nv^*-R_NDu_M^*)\|_{V_N}\\
&\leq &
(2C_A+2)\|v^*-Du_M^*)\|_{V},\\
\|v^*-D\tilde u_M\|_{D}
&\leq & 
\|v^*-Du_M^*\|_{D}+ \|Du_M^*-D\tilde u_M\|_{D}\\
&\leq & 
(2C_A+3)\|v^*-Du_M^*\|_{D}.
\end{array}
$$ 
\begin{Corollary}\RSlabel{CorConvIllPos}
Assume an analytic problem that has a MRD discretization without being
well--posed, and assume that the given data do not necessarily come from some solution $u^*$.
Then there is a uniformly stable computational strategy that 
provides trial elements that reproduce the given data at the 
quality of the {\em Trial Space Data Approximation}. This reduces the error and convergence analysis 
to an approximation problem for a data element $v^*\in V$ by a data subspace $D(U_M)$
in $V$ under the norm in $V$.\qed
\end{Corollary} 
If there is no well--posedness, there still is a {\em backward error
  analysis}. Instead of solving the problem with data $v^*$, which may
be unsolvable or ill--posed, one can come up with an element $\tilde u_M$
from the trial space which has data that are close to the given data,
and roughly as close as possible for the given trial space. For PDE solving,
this
usually means that one has  an exact solution of a PDE with perturbed boundary
data and a perturbation in the inhomogeneity of the PDE. If these perturbations
are calculated and turn out to be tolerable, the user might be satisfied with 
$\tilde u_M$. Many application papers proceed this way, unfortunately, but users should always
keep in mind that there may be very different trial elements that reproduce
the data nicely, if there is no well--posedness.
%Fixing a finite--dimensional trial space $U_M$ beforehand is a very strong form
%of {\em regularization}, but it is a crucial ingredient of Theorems
%\RSref{TheGenStab} and \RSref{TheGenConvNoise} 
%that we cannot bypass without changing the whole setting.
%We shall change the setting in the next section.
%****************************************************************
\section{Discretization in Hilbert Spaces }\RSlabel{SecDiHS}
We now assume that $U$ is a Hilbert space with inner product $(.,.)_U$
and that the data map is composed of continuous functionals
$\lambda\in\Lambda\subset U^*$ like in the beginning of section \RSref{SecIntroDD}.
The Riesz map allows a transition from functionals to functions, and thus we can 
fix a finite subset $\Lambda_N=\{\lambda_1,\ldots,\lambda_N  \}\subset \Lambda$
and consider the Riesz representers $u_1,\ldots,u_N\in U$ of these functionals.
If linear independence is assumed, we have $N$--dimensional spaces $L_N\subset U^*$ and
$U_N\subset U$  by taking the spans, and the space $V_N$ is $\R^N$ as the range
of the restriction $R_N$ with $R_NDu=(\lambda_1(u),\ldots,\lambda_N(u))^T$ which
just is the usual projection from $V:=\R^\Lambda$ to
$V_N:=\R^{\Lambda_N}=\R^N$. 
If orthonormal bases are chosen, we have the 2--norm of coefficients as
$\|u\|_U$
for all $u\in U$, but in order to comply with Section \RSref{SecIntroDD}, we
have
to take the sup--norm in the range of the data map, which is the identity if
discretized in that basis. But then the identity map is not well--posed,
due the choice of norms which is not adequate for Hilbert spaces. 
\biglf
We thus have to change the setting, taking the norms in $V_N=\R^N$ as 2-norms,
assuming $\Lambda$ to be countable and total, taking orthonormal bases, and the
restrictions as projections focusing on finite subsets of indices in the
expansions.
For the choice of $U_N$ and $V_N$ as above, we then have \eref{equMD2RNDuM}
and \eref{equMU2C}
with the constant 1.  
\biglf
This is the standard situation in  Rayleigh--Ritz--Galerkin methods.
It might be surprising that everything is perfectly well--conditioned 
here, but this is no miracle because we focused on spaces, not on bases,
and used an optimal basis for the theoretical analysis. The usual
problems with conditions of stiffness matrices etc. 
are basis--dependent, not space--dependent. 
%****************************************************************
\section{Optimal Recovery in Hilbert Spaces }\RSlabel{SecORiHS}
When starting from a finite set $\Lambda_N$ of functionals providing the data 
$\lambda_j(u^*)$ of a true solution of the analytic problem, the above choice
of a trial space as the space spanned by the representers of the functionals
is optimal under all other choices of trial spaces. This is a standard result 
in the theory of Reproducing Kernel Hilbert Spaces, but we include it here in a
general form, because of its central importance within the context
of studying all possible discretizations. 
\begin{Theorem}\RSlabel{TheOptRec}
Assume that we have a computational problem posed in a Hilbert space $U$,
and the only available data are of the form
$\lambda_1(u^*),\ldots,\lambda_N(u^*)$ for $N$ linearly independent data
functionals in $U^*$ and an unknown object $u^*\in U$. 
Then, for any linear functional $\mu\in U^*$ ,
consider all possible linear computational procedures for 
calculating good approximations of $\mu(u^*)$ using  only the above values.
Then there is a unique error--optimal strategy that works as follows:
\begin{enumerate}
 \item Use the representers $u_1,\ldots,u_N\in U$ of the functionals 
$\lambda_1,\ldots,\lambda_N\in U^*$.
\item Calculate the interpolant $\tilde u$ to $u^*$ in the span of the
  representers,
i.e. solve the system
$$
\lambda_k(u^*)=\displaystyle{\sum_{j=1}^Nc_j\lambda_k(u_j)   }
=\displaystyle{\sum_{j=1}^Nc_j(u_k,u_j)_U   } 
=\displaystyle{\sum_{j=1}^Nc_j(\lambda_k,\lambda_j)_{U^*},\;1\leq k\leq N   } 
$$ 
and define
$$
\tilde u =\displaystyle{\sum_{j=1}^Nc_ju_j}. 
$$
\item For each data functional $\mu\in U^*$, use the value 
$$
\mu(\tilde u)=\displaystyle{\sum_{j=1}^Nc_j\mu(u_j)}
$$
as an approximation to $\mu(u^*)$.
\end{enumerate} 
This approximation has minimal error under all other linear computational
procedures using the same data for calculating approximations of $\mu(u^*)$,
in the sense that the error functional has minimal norm. \qed
\end{Theorem}
This technique is independent of well--posedness and makes optimal use
of the available data, error--wise. From the previous section we conclude
that it is uniformly stable when considered in terms of spaces, not bases.
If applied to PDE solving, it is realized by {\em symmetric collocation}
\RScite{schaback:2013-2}. It can also be applied to numerical integration and
numerical differentiation, see e.g.
\RScite{davydov-schaback:2013-1,schaback:2014-1}.
\biglf
In the context of this paper, the above result shows that the quest for
good trial spaces and well--posed discretizations has a simple
solution in the Hilbert space situation. 
We shall come back to this in Section \RSref{SecWSC} when we look at the
differences between weak and strong formulations.
%****************************************************************
\section{Bases}\RSlabel{SecB}
We now assume that we have a well--posed analytic problem in the sense of Section
\RSref{SecIntroD} with an MRD discretization,
and by a proper choice of restrictions $(R_N,V_N)$ according to Theorem 
\RSref{TheGenStab}, we have uniform stability in the form of \eref{equMU2C}.   
We specialize here to the case of Theorem \RSref{TheExMRD}
where we have functionals $\lambda\in\Lambda$ 
and restrictions working via subsets $\Lambda_N\subset \Lambda$ selecting finitely
many data. 
We now choose a basis $u_1,\ldots,u_M$ of $U_M$ and take the functionals 
$\lambda_1,\ldots,\lambda_N$ from the set $\Lambda_N$.
%, and we normalize both 
% the objects and the functionals. 
Then we
consider the discretized system 
\bql{eqaRufDisc}
\sum_{j=1}^Ma_j\lambda_k(u_j)\approx f_{\lambda_k}=\lambda_k(u^*),\,1\leq k\leq N
\eq
that we solve approximatively by residual minimization like in Section 
\RSref{SecIntroSDP}.
Clearly, a bad choice of bases will spoil stability, but we want to 
study this effect in detail.
We quantify the stability of the object basis by norm equivalence
$$
\begin{array}{rclcll}
c_M\|u_a\|_{WP} &\leq & \|a\|_M &\leq & C_M\|u_a\|_{WP} & \fa a\in\R^M,\\
%c_N\|\lambda_b\|_{U^*} &\leq & \|b\|_N &\leq & C_N \|\lambda_b\|_{U^*}& \fa b\in\R^N.
\end{array}
$$
with an unspecified norm $\|.\|_M$ on $\R^M$ that is used in computation.
% and $\|.\|_N$ on $\R^N$. 
With the $N\times M$ matrix $A=(\lambda_k(u_j))_{1\leq j\leq M,\,1\leq k\leq N}$. 
and the basis representation
$$ 
\begin{array}{rcl}
u_a&:=&\displaystyle{\sum_{j=1}^Ma_ju_j   } \\
%\lambda_b&:=&\displaystyle{\sum_{k=1}^Nb_k\lambda_k   } 
\end{array}
$$
with coefficient vectors $a\in\R^M$, % and $b\in\R^N$, 
we see that $Aa=R_NDu_a$ holds. and get
$$
\begin{array}{rcl}
\|a\|_M
&\leq &
C_M\|u_a\|_{WP} \\
&\leq &
CC_M\|u_a\|_D \\
&\leq&
2CC_M\|R_NDu_a\|_{V_N} \\
&=&
2CC_M\|Aa\|_{V_N} \\
\end{array}
$$
by \eref{equWPCDu}, and \eref{equMD2RNDuM}.
\begin{Theorem}\RSlabel{TheMatStab}
Under the above assumptions, the system \eref{eqaRufDisc} 
has the stability property
$$
\|a\|_M\leq 2CC_M\|Aa\|_{V_N} \fa a\in\R^M.\qed
$$
\end{Theorem}
In Section \RSref{SecIntroSDP}, we minimized $\|R_N(Du^*-Du_M)\|_{V_N}$ over all
$u_M\in U_M$. After introducing a basis in $U_M$, this is the same as minimization
of $\|f-Aa\|_{V_N}$ with $f:=R_NDu^*=(\lambda_1(u^*),\ldots,\lambda_N(u^*))^T\in\R^N$
over all $a\in\R^M$. We are satisfied with a vector $\tilde a\in\R^M$ such that
\bql{eqfAaVN}
\|f-A\tilde a\|_{V_N}\leq \min_{a\in\R^M}\|f-Aa\|_{V_N}\leq C_A\|f-Aa^*\|_{V_N},
\eq
where $a^*$ is a good coefficient vector for the direct approximation of
the true solution $u^*$ by elements of the trial space $U_M$.
We use $a^*$ in the way we used $u_M^*$ in Section \RSref{SecIntroSDP}
as a competitor that may come from some special approximation technique.
Then we form the elements $\tilde u_M:=u_{\tilde a},\;u_M^*=u_{a^*}\in
U_M$ 
and see that \eref{equuM} is satisfied. 
\biglf
This lets us arrive at
Theorem \RSref{TheGenConv}, implying that the convergence rate is the same as
the rate for the {\em Trial Space Data Approximation},
but this does not yield error
bounds
in terms of coefficients. However, we can proceed by
$$
\begin{array}{rcl}
\dfrac{1}{2CC_M}\|a^*-\tilde a\|_M
&\leq & 
\|A(a^*-\tilde a)\|_{V_N}\\
&\leq & 
\|Aa^*-f\|_{V_N}+\|f-A\tilde a\|_{V_N}\\
&\leq & 
(1+C_A)\|Aa^*-f\|_{V_N}\\
\end{array} 
$$
and get
$$
\|a^*-\tilde a\|_M\leq (1+C_A)(2CC_M)\|Aa^*-f\|_{V_N}.
$$
The norm in $V_N$ must be chosen to comply with Section \RSref{SecIntroDD},
and this works for the discrete sup norm. But if users do not want to 
minimize $f-Aa$ in the sup norm, an additional norm equivalence
comes into play,  now on $V_N$, and this will often depend on $\dim V_N$.
In detail, norm equivalence in $V_N$ is assumed as 
$$
c_N\|R_Nv\|_{V_N} \leq  \|R_Nv\|_N \leq  C_N \|R_Nv\|_{V_N} \fa v \in V.
$$
and minimization in the new norm $\|.\|_N$ will replace \eref{eqfAaVN} 
by
$$
\|f-A\tilde a\|_{N}\leq \min_{a\in\R^M}\|f-Aa\|_{N}\leq C_A\|f-Aa^*\|_{N},
$$
and our above argumentation now yields
$$
\|a^*-\tilde a\|_M\leq \dfrac{(1+C_A)(2CC_M)}{c_N}\|Aa^*-f\|_{N}.
$$
If bases are chosen badly, the quotient $C_M/c_N$ can be extremely large 
and will spoil the uniformity that we had so far.
\biglf
Users can check their stiffness matrices $A$ computationally for stability,
but Theorem \RSref{TheMatStab} indicates that there may be a strong influence 
due to a bad choice of the trial basis. Even a calculation of a
Singular Value Decomposition will not be completely basis--independent, since it
only eliminates {\em orthogonal} basis transformations in the domain and range of $A$.
%****************************************************************
\section{Nodal Bases}\RSlabel{SecNod}
In {\em meshless methods}, it is customary to write everything 
``entirely in terms of nodes'' \RScite{belytschko-et-al:1996-1},
which means that the functions $u_M$ in the trial space $U_M$ are parametrized by
their values at certain {\em nodes} $x_1,\ldots,x_M$, i.e.
$$
u_M(x)=\displaystyle{\sum_{j=1}^M s_j(x)u_M(x_j) \fa u_M\in U_M   } 
$$ 
with {\em shape functions} $s_j$ that are usually localized around $x_j$
and have the Lagrange property $s_j(x_k)=\delta_{jk},\,1\leq j,k\leq M$.
We prefer the term {\em nodal basis}, because there is nothing meshless 
in the above representation, and the standard finite elements, which nobody
would call meshless, are nodal as well in the above sense. 
Many application papers report experimentally that these bases have favorable stability
properties, and we shall now show why. 
\begin{Theorem}\RSlabel{TheGenNodalStab}
Assume a well--posed problem in the sense of \eref{equWPCDu}, 
where $U$ is a space of functions on some domain $\Omega$. Furthermore, 
assume that the point evaluation functionals $\delta_{x}$ are uniformly bounded 
by $\gamma>0$
in the norm $\|.\|_{WP}$. Finally, assume that the data space $V$ and the
restrictions $V_N$ are normed via supremum norms, as mentioned in 
Section
\RSref{SecIntroDD} and Theorem \RSref{TheExMRD} as a special case. 
Then for each trial space $U_M\subset U$
with a nodal basis $s_1,\ldots,s_M$ using nodes $x_1,\ldots,x_M\in \Omega$ one can find a finite
set
of functionals $\lambda_1,\ldots,\lambda_N$ such that the $N\times M$ {\em stiffness matrix}
$A$ with entries $\lambda_j(s_k)$ has the uniform stability property
\bql{eqa2gCAa}
\|a\|_{\infty}\leq 2\gamma C \|Aa\|_{\infty} \fa a\in \R^M.
\eq
\end{Theorem} 
{\bf Proof}:  We apply Theorem \RSref{TheGenStab}.  Then
$$
|u_M(x_j)|\leq \gamma \|u_M\|_{WP}\leq 2\gamma C\|R_ND(u_M)\|_{V_N}
= 2\gamma C \displaystyle{\max_{\lambda_k\in \Lambda_N}\left|\sum_{j=1}^M \lambda_k(s_j)u_M(x_j)\right|  } 
$$
and if we denote the vector of nodal values by $u_X\in \R^M$, we see that
$$
\|u_X\|_\infty\leq 2\gamma C \|Au_X\|_\infty
$$ 
with the {\em stiffness matrix} $A$ with entries $\lambda_k(s_j)$. \qed
\biglf
This means that all trial spaces with nodal bases can be uniformly
stabilized by taking good and large selections of test functionals. Furthermore, Section
\RSref{SecIntroSDP}
provides convergence proofs and convergence rates for such techniques. 
\biglf
We now go closer to what user would do. In the notation of Section
\RSref{SecB} they would invoke a least--squares solver minimizing
$\|f-Aa\|_2$ instead of minimizing $\|R_ND(u^*-u_M)\|_{V_N}$ which is
$\|f-Aa\|_\infty$ in terms of linear algebra. In the notation of Section
\RSref{SecB},
we then have $\|.\|_M=\|.\|_{\infty,\R^M}$, $\|.\|_N=\|.\|_{\ell_2,\R^N}$,
$C_M=\gamma$, $\|.\|_{V_N}=\|.\|_{\infty,\R^N}$, $c_N=1$, leading to 
$$
\|a^*-\tilde a\|_\infty\leq 2\gamma C(1+C_A)\|Aa^*-f\|_{2}
$$
for any reference approximation $a^*$. In this case, we may take $a^*$ as the
vector of nodal values of the true solution, and then 
\bql{eqMaxNodeErr}
\displaystyle{\max_{1\leq j\leq M} |u^*(x_j)-\tilde u_M(x_j)|
\leq 2\gamma C(1+C_A) \max_{1\leq k\leq N}|\lambda_k(u^*-s^*)|}
\eq
where $s^*$ is the trial function with the nodal values of the true solution.
\begin{Corollary}\RSlabel{CorNodalConv}
Assume a well--posed analytic problem with an MRD discretization, and assume
that a trial space is parametrized by a nodal basis. Then the error of a
computational procedure as in Section \RSref{SecIntroSDP}, evaluated  on the nodes, 
is pointwise bounded by the error of the {\em Trial Space Data Approximation},
i.e. the approximation of the data of the true solution by
the data of trial elements, measured in the data norm.\qed 
\end{Corollary}
This means that using a nodal basis transfers the results of Theorem
\RSref{TheConv2} directly to a convergence on the nodes. This is a very useful
result for many meshless methods using nodal bases, e.g. when applying
Moving Least Squares techniques.
%****************************************************************
%****************************************************************
%****************************************************************
\section{Examples}\RSlabel{SecEx}
%****************************************************************
\subsection{Interpolation}\RSlabel{SubSecInt}
For illustration, we start with the rather simple case of
recovering a function $u$ on some compact domain $\Omega\subset \R^d$ from data 
of $u$ that do not involve derivatives. A {\em strong} formulation 
takes  $\Lambda=\{\delta_x\;:\;x\in \Omega\}$ on a space $U$
on which these functionals are continuous, e.g. $U=C(\Omega)$ under the sup
norm. A {\em weak} formulation uses different data, e.g. functionals
$$
\lambda_v(u):=(u,v)_{L_2(\Omega)} \fa u,v
\in U:=L_2(\Omega)
$$
and 
$$
\Lambda:=\{\lambda_v\;:\;\|v\|_{L_2(\Omega)}=\|\lambda_v\|_{L_2^*(\Omega)}=1  \}.
$$
The strong case takes $V:=C(\Omega)=U$ under the sup norm, while the weak case 
uses $V:=L_2(\Omega)=U$ under the $L_2$ norm. In both cases, the data map is the
identity, and we have well--posedness in the norms $\|.\|_{WP}=\|.\|_U$ in both
cases,
but the norms differ.
\biglf
The restrictions can
work by selection of finitely many functionals in both cases, and all axioms of
Section 
\RSref{SecIntroDD} are satisfied. 
\biglf
We now fix an arbitrary finite--dimensional trial space $U_M\subset U$,
and Theorem \RSref{TheGenStab} yields that there is a restriction that makes the
linear system \eref{eqRNDUM} uniformly stable in the sense \eref{equMU2C} with $C=1$.
\biglf
The computational procedures of Section \RSref{SecIntroSDP} can use a comparison
trial object $u_M^*$ that is the best approximation to the true solution $u^*$
in the norm $\|.\|_U=\|.\|_{WP}=\|.\|_D$, and Theorem \RSref{TheGenConv} then
shows that the computational solution $\tilde u_M$ has the same convergence rate
as the 
best approximation $u_M^*$.  
\biglf
The computational solution $\tilde u_M$ is obtained via \eref{equuM} from a stably discretized
linear system, and we assume that we perform inexact minimization of 
$\|R_N(u^*-u_M)\|_{V_N}$. 
\biglf
In the strong case, this is best linear discrete Chebyshev
approximation on sufficiently many points, i.e. a linear optimization problem. 
In the univariate case with $U_M$ being a trial space of polynomials of degree
$M-1$ on an interval $I$, a discretization on $N\geq M $ test points forming a set $P_N$
will always have a stability  inequality
$$
\|u\|_{\infty,I}\leq C(M,N)\|u\|_{\infty, P_N}
$$
of the form \eref{equMDDuMV}, but the stability constant varies.
For $M=N$ equidistant points, the constant $C(M,M)$ grows exponentially with
$M$, and for Chebyshev--distributed test points it still grows like $\log M$.  
Uniform stability holds for $N={\calO}(M^2)$ equidistant points, as follows from
a standard argument going back to the notion of 
norming sets \RScite{jetter-et-al:1999-1} 
and using Markov's inequality \cite[Ch. 3.3]{wendland:2005-1}.
Theorem \RSref{TheGenStab} only proves {\em existence} of a uniformly stable
discretization,
but this example shows that there may be a considerable amount of oversampling
or overtesting behind the scene.
\biglf
If stability is uniform, nodal bases written in terms of values at $M$ nodes
$x_j$ will trivially lead to $|u(x_j)|\leq  \|u\|_{\infty,I}\leq C\|u\|_{\infty, P_N}\fa u\in U_M$, which is 
\eref{eqa2gCAa}. 
\biglf
The weak case discretizes by $N$ well--chosen normalized test functionals $\lambda_{v_j}$ with
normalized Riesz representers $v_j\in L_2(\Omega)$, and the quantity
$\|R_N(u^*-u_M)\|_{V_N}$ to be minimized is 
$$
\displaystyle{ \max_{1\leq j\leq N}|(v_j,u^*-u_M)_{L_2(\Omega)}|   }.
$$ 
Our theory shows that the test functionals can be chosen to render uniform
stability,
but there is a trivial standard choice via the $M$ functionals represented
by an orthonormal basis
$v_1,\ldots,v_M$ of $U_M$. Then the above minimization produces the best
approximation $u_M^*$ to $u^*$ from $U_M$ without any oversampling.
\biglf
In both cases, Theorem \RSref{TheConv2} is applicable, and we see that we the
$L_\infty$
or $L_2$ convergence rates of the non--discrete best approximations carry over
to
the discrete approximations.
\biglf
To compare the difference of convergence rates between weak and strong formulations for a
given fixed trial space $U_M$, we see immediately that the $L_2$
convergence rate  is never worse than the $L_\infty$ rate, but it is taken in a
weaker norm. If users insist on the best possible convergence rate in $L_2$,
they should take a weak form, at the expense of a sufficiently good numerical
integration. But the $L_\infty$ error of their solution will clearly
not have a better $L_\infty$ convergence rate than the strong solution.
\biglf
Both computational approximations, weak or strong, converge like the best
approximations in the respective data norm, and this is a fair deal.
Convergences can be spectral  in certain cases, e.g.
in case of univariate functions on an interval $I$ that have a complex extension that is analytic in a 
region of the complex plane containing $I$ in its interior. 
This shows how the theory applies to spectral convergence situations without change.
\biglf
But, of course, there is the extreme case where the solution is only in $L_2$
and not in $C(\Omega)$. Then the strong technique is undefined. But then the weak
technique is forced to have the weak data given directly, without numerical
integration,
because the latter is as unfeasible as the strong technique.  
\biglf
All other examples will show a very similar behavior, differing only in their
data maps.
%****************************************************************
\subsection{Standard Homogeneous Weak Poisson Problem}\RSlabel{SecSDP}
We fix a bounded Lipschitz domain $\Omega\subset\R^d$ and consider the weak
Dirichlet problem $-\Delta u=f$ with homogeneous boundary conditions.
This works on the Hilbert space $U:=H_0^1(\Omega)$ with the inner product
$$
(u,v)_1:=\int_\Omega \nabla^Tu(x) \nabla v(x)dx \fa u,v\in U,
$$
and the standard (global) weak formulation asks for a function $u\in U$ with
$$
(u,v)_1=(f,v)_{L_2(\Omega)}  \fa v\in U=H_0^1(\Omega).
$$
In the sense of this paper, the functionals are 
$$
\lambda_v\;:\;u\mapsto \lambda_v(u)=(u,v)_1 \fa u,v\in U 
$$
and the problem takes the form \eref{eqruf} with 
$$
\begin{array}{rcl}
\Lambda &:=& \{ \lambda_v\;:\;v\in U,\;\|v\|_U=1\}\subset U^*\\
f_{\lambda_v} &=& (u^*,v)_1=(f,v)_{L_2(\Omega)}\fa \lambda_v\in \Lambda
\end{array}
$$
where $u^*\in U$ is the true solution.
\biglf
To check the well--posedness in the sense of section \RSref{SecIntroD},
we get
$$
\begin{array}{rcl}
\displaystyle{  \|Du\|_V=\sup_{\lambda_v\in \Lambda}|\lambda_v(u)|
=\sup_{v\in U,\|v\|_U=1}|(u,v)_U|=\|u\|_U=:\|u\|_{WP}}
\end{array}
$$
proving well--posedness, and the data map $D$ is an isometry.
\biglf
We now consider fairly arbitrary trial spaces $U_M\subset U=H_0^1(\Omega)$
to allow standard or extended or generalized finite elements, or even
certain spectral methods of Galerkin type. Theorem \RSref{TheConv2} is
applicable, and we see that we get the convergence rate of approximations to the
true solution in $U=H_0^1(\Omega)$. This is well--known from finite elements,
but it holds in general, provided that MRD testing is done. 
It applies to Petrov--Galerkin methods and 
spectral techniques of Galerkin type. The rate 
mainly depends on the smoothness of the solution and on the
trial space chosen.
\biglf
For the standard finite--element situation with piecewise linear elements,
this yields ${\calO}(h)$ convergence in $H_0^1(\Omega)$, as usual for that regularity. To reach
${\calO}(h^2)$ convergence in $L_2(\Omega)$ under $H^2(\Omega)$ regularity, the Aubin--Nitsche
trick is an add--on that is not covered by our theory. But it follows from the
fact that the best approximation to $u^*$ in $H_0^1(\Omega)$ automatically has
${\calO}(h^2)$ convergence in $L_2$ under $H^2(\Omega)$ regularity. This is independent
of PDE solving, it is a property of Approximation Theory. 
%****************************************************************
\subsection{Collocation Methods}\RSlabel{SubSecCM}
We now want to focus on the general statement
\begin{itemize}
 \item[] All linear PDE or ODE problems can be 
numerically solved by collocation in sufficiently many points in such a way  
that  
the convergence rate in the well--posedness norm is at least the rate of the 
{\em trial space
data approximation}.
\end{itemize}
This, of course, includes pseudospectral methods. 
But we have to add more details to show how  it 
follows from Theorem \RSref{TheConv2}.  We only have to show that collocation
is an MRD discretization and pick a suitable form of well--posedness.
\biglf
The space $U$ should be a normed linear space of functions on a domain
$\Omega$ with boundary $\Gamma$, for instance a Sobolev space. To keep things
simple,
we assume that 
the analytic problem is posed in strong form by evaluating 
a linear elliptic second--order differential operator $L$ on points of the domain and 
a linear boundary operator $B$ on the boundary, i.e.
\bql{eqLBProb}
\begin{array}{rcll}
Lu(x) &=& f(x), &\fa x\in \Omega\\
Bu(y) &=& g(y), &\fa y\in \Gamma\\
\end{array}
\eq
where $f$ and $g$ are given functions on $\Omega$ and $\Gamma$.
Introducing continuous functionals $\lambda_x(u):=Lu(x)=\delta_x \circ L$ and 
$\mu_y(u):=Bu(y)=\delta_y \circ B$ on $U$ one gets a problem of the form 
\eref{eqruf} with
\bql{eqLamlam}
\Lambda:=\{\lambda_x\;:\;x\in \Omega\}\cup\{\mu_y\;:\;y\in \Gamma\} 
\eq
and it should be clear that one can allow more than two 
operators, and combinations of different boundary conditions.
\biglf 
From here there are different ways to proceed towards well--posedness,
but we can ignore well--posedness for a moment. 
We normalize all functionals as elements of $U^*$ and pose the problem
in the form  \eref{eqruf} with $f_\lambda:=\lambda(u^*) \fa \lambda\in \Lambda$.
Then $|\lambda(u)|\leq \|u\|_U \fa u\in U,\;\lambda\in \Lambda$,
and there is no problem to define the space $V$ and the restrictions via 
taking suprema. We can apply Theorems \RSref{TheGenStab} and \RSref{TheGenConv}
without assuming well--posedness, and we see that we can work on any trial space
$U_M$, but all results only hold in the data norm. Comparing with any existing good
approximation $u_M^*$ to $u^*$ from $U_M$, we get some $\tilde u_M\in U_M$ 
by a discrete computational method such that
$$
\displaystyle{ \sup_{\lambda\in \Lambda}|\lambda(\tilde u_M-u^*)|
\leq (2C_A+3)\sup_{\lambda\in \Lambda}|\lambda(u_M^*-u^*)| \leq (2C_A+3)\|u_M^*-u^*\|_U  }
$$
due to normalization of the functionals, and if we use some   $u_M^*$ 
with small $\|u_M^*-u^*\|_U$, we get the above statement. The backward error
analysis of Section \RSref{SecIntroND} will be applicable here.
\biglf
In the above setting, the most natural well--posedness condition would be
of the form
\bql{equUCLuBu} 
\|u\|_{WP}\leq C \max(\|Lu\|_{\infty,\Omega},\|Bu\|_{\infty,\Gamma})
\eq
for a suitable norm $\|.\|_{WP}$ on $U$. This holds for
$U:=C^2(\Omega)\cap C(\overline \Omega)$ with the sup norm in $U$ 
\cite[(2.3), p. 14]{braess:2001-1} for uniformly elliptic operators $L$ and
Dirichlet boundary data. This implies by Theorem \RSref{TheConv2} that for such problems the 
convergence rate of the {\em Trial Space Data Approximation} carries over 
to the convergence rate of collocation in the
sup norm.
\biglf
As a special example we consider unsymmetric collocation \RScite{kansa:1986-1} by translates 
of the kernel $K$ of a Hilbert space $U$, applied to a Dirichlet problem of the form
\eref{eqLBProb} on a domain $\Omega\subset\R^d$. The trial
space $U_M$ is spanned by kernel translates $v_j:=K(\cdot,x_j),\;1\leq j\leq M$ for nodes
$x_1,\ldots,x_M\in \Omega$, but this is not a stable basis. A {\em nodal basis} 
in the sense of Section \RSref{SecNod} 
consists of the Lagrange basis $u_1,\ldots,u_M$ spanning the same trial space.
Collocation is done via the functionals defined for \eref{eqLamlam}, and to make
them continuous we can take a space like $U:=H^m(\Omega)$ with some $m>2+d/2$. We have
well--posedness in the sense of \eref{equUCLuBu} in the sup--norm. 
\begin{Theorem}\RSlabel{TheAsyColl}
Unsymmetric collocation in the sense of E. Kansa \RScite{kansa:1986-1}
has the property that for each possible trial space spanned by kernel translates 
there is a selection of test functionals such that the stiffness matrix,
when written in terms of the nodal basis, 
has a uniform stability property \eref{eqa2gCAa}. If solved by residual
minimization along the lines of Section \RSref{SecIntroSDP}, error bounds follow from
Corollary \RSref{CorConv} or \eref{eqMaxNodeErr}. Convergence rates 
in the sup norm are obtained from the rate of convergence of second derivatives 
in the sup norm of interpolants of the true solution by the trial space.  \qed
\end{Theorem}  
This provides many explicit convergence rates via standard results on 
interpolation by translates of kernels \cite[Chapter 11]{wendland:2005-1}. For
instance, the convergence for the Whittle--Mat\'ern kernel reproducing
$H^m(\R^d)$ for $m>2+d/2$ is like ${\calO}(h^{m-2-d/2})$ in terms of the {\em fill distance} 
$h:=\sup_{y\in \Omega}\min_{1\leq j\leq M}\|y-x\|_2$, while the convergence is
exponential 
for kernels like the Gaussian or multiquadrics.   
\biglf
The functionals in \eref{eqLamlam} are a mixture of two kinds, but 
Theorem \RSref{TheGenStab} and Lemma \RSref{LemDatStab}  do not say 
how to achieve a uniformly stable 
balance between testing $B$ on the boundary and testing $L$ in the interior.
Future work should address this problem, and Section \RSref{SubSecInt} 
suggests that there might be quite some overtesting needed for uniform
stability.
Square collocation systems can even be singular \RScite{hon-schaback:2001-1},
such that overtesting is necessary in general.
\biglf
All of this readily generalizes to plenty of other linear well--posed PDE
problems,
and readers can use the tools of this paper to assemble what they need.
Note that unsymmetric collocation is a pseudospectral method 
in the sense of the literature 
(see e.g. \RScite{fornberg-sloan:1994-1,fornberg:1996-1,canuto-et-al:2007-1}) on
spectral methods,
and this paper provides a general way to assess convergence of pseudospectral
methods.
Since we write the analytic and the computational problems in terms of arbitrary
functionals, this approach also covers spectral methods in Tau form.
%****************************************************************
\subsection{Weak Dirichlet Problems}\RSlabel{SecWDP}
The standard finite element procedures for solving Dirichlet problems for the
Laplace operator on bounded domains $\Omega\subset\R^d$ use 
strong data on the boundary and weak data in the interior. The data functionals
are
$$
\begin{array}{rcl}
\Lambda_1
&:=& \{\lambda_v\;:\;u\mapsto (\nabla u, \nabla v)_{L_2(\Omega)} \fa v\in
H_0^1(\Omega),\;\|\nabla v\|_{L_2(\Omega)}=1 \},\\
\Lambda_2&:=& \{\delta_y\;:\;y\in \Gamma:=\partial \Omega\},\\
\Lambda&:=& \Lambda_1\cup\Lambda_2.
\end{array} 
$$ 
This leads to the data norm
$$
\|u\|_D=\max(\|u\|_{\infty,\Gamma}, \|\nabla u\|_{L_2(\Omega)})
$$
if we take the sup over all functionals as in Section \RSref{SecIntroDD}.
It is well-defined on the space $U:=H^1(\Omega)\cap C(\overline\Omega)$.
Using the Poincar\'e inequality and the Maximum Principle \RScite{jost:2002-1}
after splitting $u$ into a harmonic part with boundary conditions and
a function in $H_0^1(\Omega)$ satisfying the differential equation, we get 
a well--posedness inequality
$$
\|u\|_{L_2(\Omega)}\leq C\|u\|_D =C\max(\|u\|_{\infty,\Gamma}, \|\nabla u\|_{L_2(\Omega)})\fa u\in U.
$$
Note that the Sobolev inequality forbids to use the sup norm
on the left--hand side for space dimension $d>1$.
\biglf
Whatever the chosen trial spaces in $U$ are, Theorem \RSref{TheConv2} 
shows that the convergence rate in $L_2$ of 
uniformly stabilized computational methods will be the convergence rate of 
the {\em Trial Space Data Approximation}, i.e. with respect to
$\|u\|_{\infty,\Gamma}$ and $\|\nabla u\|_{L_2(\Omega)}$.
If the trial space is spanned by translates of the Whittle--Mat\'ern kernel reproducing
$H^m(\R^d)$ for $m>1+d/2$, the rate is ${\calO}(h^{m-1-d/2})$ in terms of the fill distance
of the trial nodes \RScite{wendland:2005-1}. 
\biglf
For standard finite elements, the above approach yields ${\calO}(h)$ convergence in $L_2$.
This is without the Aubin--Nitsche trick, and it does not use $H^2$ regularity.
\biglf
The Aubin--Nitsche trick has nothing to do with finite elements and weak problems.
It is a feature of Approximation Theory,
doubling a convergence rate for certain nested approximations
in Hilbert spaces under additional regularity
assumptions. This is well--known from splines 
\cite[5.10]{ahlberg-et-al:1967-1}
and kernel--based methods \RScite{schaback:1999-1}.
In the context of this paper, one considers the best approximation to the true
solution in $H_0^1(\Omega)$, and it will automatically yield ${\calO}(h^2)$
convergence under $H^2$ regularity, but only for zero boundary conditions. 
%****************************************************************
\subsection{Weak--Strong Comparison}\RSlabel{SecWSC}
If we compare with what we had in the strong case, the situation 
for fixed trial spaces is roughly as follows:
\begin{enumerate}
 \item The weak case has $L_2$ convergence at the convergence rate for first
   derivatives,
\item the strong case has $L_\infty$ convergence at the convergence rate for second derivatives.
\end{enumerate}
This usually yields a slightly better rate for the weak case,
as we saw when comparing ${\calO}(h^{m-1-d/2})$ with ${\calO}(h^{m-2-d/2})$ 
for the trial space spanned by translates of the
Whittle--Mat\'ern kernel.
On the downside, weak methods usually have to perform numerical integration
at an accuracy that complies with the convergence rate, 
and they converge in a weaker norm.
\biglf
If one fixes the available finite data and then looks for an error--optimal 
solution in a fixed Reproducing Kernel Hilbert Space, the above discussion about 
differences between strong and weak methods becomes obsolete. The optimal
solution is always the one described in Section \RSref{SecORiHS}, and it is
furnished by symmetric collocation \RScite{schaback:2013-2}. Since it allows
arbitrary evaluation functionals $\mu$ in Theorem \RSref{TheOptRec}, it is
pointwise and $L_\infty$--optimal by taking functionals $\mu=\delta_x$, and
$L_2$--optimal by taking all functionals $\mu=\lambda_v=(.,v)_{L_2(\Omega)}$. 
%****************************************************************
\subsection{MLPG}\RSlabel{SecMLPG}
%****************************************************************
We stay with the Dirichlet problem for the Laplacian, for simplicity, 
and describe the standard variation of the 
{\em Meshless Local Petrov--Galerkin} \RScite{atluri-zhu:1998-1,atluri:2005-1} method. 
The difference to the standard weak formulation is that the integrals are
localized and the boundary integrals are kept. This means that
on small subdomains $\Omega_h\subset\Omega$ 
with boundaries $\Gamma_h\subset \overline \Omega$ the strong equation
$-\Delta u=f$ is integrated against a test function $v_h$ to define functionals
of the form 
$$
u\mapsto \lambda_{\Omega_h,v_h}(u)=\dfrac{-1}{vol(\Omega_h)}\int_{\Omega_h}v_h\cdot \Delta u
$$
that are continuous on $U:=C^2(\Omega)\cap C(\overline{\Omega})$, and the
problem \eref{eqLBProb} takes the form \eref{eqruf} via
\bql{eqlamdg}
\begin{array}{rcll}
\lambda_{\Omega_h,v_h}(u)
&=& 
\dfrac{1}{vol(\Omega_h)}\int_{\Omega_h}v_h\cdot f \;\; &\fa \Omega_h\subset \Omega,
\;v_h\in C(\Omega_h)\\
\delta_y(u) &=& g(y) &\fa y\in \Gamma
\end{array}
\eq
for given continuous functions $f$ on $\Omega$ and $g$ on $\Gamma$. One can
restrict the domains $\Omega_h$ and the test functions $v_h$ further, and allow
other ways of handling the boundary conditions. Furthermore, the above
functionals
are usually transformed by integration by parts before they are implemented, 
but we deal with this later.
\biglf
The goal is to prove some form of well--posedness for the analytic problem,
and this seems to be missing completely in the rich literature on the MLPG
method. On the space $U=C^2(\Omega)\cap C(\overline{\Omega})$ we know that
\eref{equUCLuBu} holds for $\|.\|_{WP}=\|.\|_{\infty,\overline \Omega}$, and we assert 
\bql{equCluu} 
\displaystyle{ \|u\|_{\infty,\overline\Omega}\leq
  C\left(\sup_{\Omega_h,v_h}
|\lambda_{\Omega_h,v_h}(u)|  +\|u\|_{\infty,\Gamma}\right)}\fa u\in U.
\eq
But this follows from \eref{equUCLuBu} by setting $f:=\Delta u$ in
\begin{Lemma}\RSlabel{LemFCont}
For each continuous function $f$ on some compact domain $\Omega\subset\R^d$ 
the norms
$$
\|f\|_{\infty,\Omega}=\|f\|_I:=\displaystyle{\sup_{\Omega_h\subset\Omega}\dfrac{1}{vol(\Omega_h)}
\left|\int_{\Omega_h}f  \right|   } 
$$
coincide, where  the diameter of the admissible sets $\Omega_h$ can be bounded above by some
arbitrary $r>0$, if required. One can also restrict the subdomains
$\Omega_h$ to
balls or cubes intersected with $\Omega$. 
\end{Lemma} 
{\bf Proof}: Clearly $\|f\|_I\leq \|f\|_{\infty,\Omega}$ holds. To prove
$\|f\|_I\geq \|f\|_{\infty,\Omega}$,
assume $f\neq 0$ with $\|f\|_{\infty,\Omega}=f(\tilde x)>0$ for some $\tilde
x\in\Omega$.   Then pick an arbitrary $\epsilon < f(\tilde x)/2$ and 
an arbitrary $r>0$ and choose
$\Omega_h$ to be a subdomain of the set of points $x\in\Omega$ with 
$$
0<(1-\epsilon)f(\tilde x)\leq f(x)\leq  f(\tilde x),\,\|x-\tilde x\|_2\leq r.
$$
For instance, one can take the intersection of sufficiently small balls or cubes around
$\tilde x$ with the domain $\Omega$, or if $\tilde x$ is on the boundary, one
may move slightly into the interior and ensure $\Omega_h$ to be in the interior
of the domain.
Then
$$
(1-\epsilon)\|f\|_{\infty,\Omega}= (1-\epsilon)f(\tilde x)\leq \dfrac{1}{vol(\Omega_h)}
\int_{\Omega_h}f \leq  f(\tilde x)=\|f\|_{\infty,\Omega}.\qed
$$
Note that this proves well--posedness only on  $U=C^2(\Omega)\cap
C(\overline{\Omega})$,
not on a larger space, but for all possible test functions and domain shapes
and sizes. The boundary conditions can be rephrased by weak functionals taking means,
using 
Lemma \RSref{LemFCont} again, now setting $f:=g$ and working on the boundary.
\biglf
Any a--priori renormalization of {\em all} available functionals will possibly
spoil this
argument. But as soon as finitely many functionals are selected for computation, one can
renormalize for the computational procedure. 
\biglf
If integration by parts is applied to the functionals, they change their
computational form without changing their value, and this is used in the known variations of the 
MLPG technique. For instance, 
$$
\begin{array}{rcl}
\lambda_{\Omega_h,v_h}(u)
&=&
\dfrac{1}{vol(\Omega_h)}\int_{\Omega_h}\nabla^Tv_h\cdot
\nabla u-\dfrac{1}{vol(\Omega_h)}\int_{\Gamma_h}v_h
\dfrac{\partial u}{\partial n}\\
&=&
\dfrac{-1}{vol(\Omega_h)}\int_{\Omega_h}\Delta v_h\cdot
u
+\dfrac{1}{vol(\Omega_h)}\int_{\Gamma_h}u
\dfrac{\partial v_h}{\partial n}
-\dfrac{1}{vol(\Omega_h)}\int_{\Gamma_h}v_h
\dfrac{\partial u}{\partial n}\\
\end{array} 
$$ 
are two ways to rewrite the functionals on different domains with different
admissible test functions. The basic well--posedness on $U=C^2(\Omega)\cap
C(\overline{\Omega})$ will stay as is, because the sup of all these functionals 
will be bounded above by $\|\Delta u\|_{\infty,\Omega}$, as long  as there are
no other upper bounds proven.
\biglf
The method called MLPG5 uses constant test functions like in Lemma
\RSref{LemFCont}. Then the functionals take the extremely simple 
form
$$
\lambda_{\Omega_h,v_h}(u) =-\dfrac{1}{vol(\Omega_h)}\int_{\Gamma_h}
\dfrac{\partial u}{\partial n},
$$
i.e. they are only integrals of the normal derivative on subdomain boundaries.
Nevertheless, Lemma \RSref{LemFCont} holds, and there is well--posedness
in the sup norm on $U=C^2(\Omega)\cap
C(\overline{\Omega})$. It is an open problem to prove other well--posedness
inequalities
after fixing a special form of the functionals. The above technique via Lemma \RSref{LemFCont}
always goes back to \eref{equUCLuBu}, whatever the form of the functionals is
after integration by parts.
Therefore the convergence theory for given trial spaces will be the same as for
the strong collocation methods in Section \RSref{SubSecCM}. 
\begin{Theorem}\RSlabel{TheMLPG}
If the Meshless Local Petrov--Galerkin
method is carried out 
\begin{enumerate}
 \item
for a well--posed second--order elliptic problem, 
\item using  sufficiently many well--chosen test functionals \eref{eqlamdg}
along the lines of Theorem  \RSref{TheGenStab},
\item and applying a residual minimization algorithm  as in Section
\RSref{SecIntroSDP} for   solving the overdetermined system approximatively, 
\end{enumerate}  
the algorithm is convergent with uniform stability, and the convergence rate 
in the sup norm is
the rate of the {\em Trial Space Data Approximation}. 
This rate is at least as good as for strong
collocation using the same trial spaces. \qed
\end{Theorem} 
Depending on the PDE problem, the smoothness of the true solution,
and the trial space chosen, this yields various  
convergence results, up to spectral convergence. In most applications, the trial
functions are shape functions provided by Moving Least Squares, and raising the
degree  of the local polynomials will increase the convergence rate
appropriately
\RScite{levin:1998-1,wendland:2000-1,armentano:2001-1,armentano-duran:2001-1}. 
Readers are encouraged to apply the framework of this paper to 
derive special convergence results for various trial spaces and different
variations of the MLPG technique. In particular, an extension to elasticity
problems should be quite useful.
\biglf
But the methods of this paper always assume the functionals to be given exactly,
not approximately. Only their values can be noisy, as in Section
\RSref{SecIntroND}.
This excludes various interesting applications, namely the Direct Meshless Local
Petrov Galerkin (DMLPG) technique \RScite{mirzaei-schaback:2011-1}
and localized kernel--based methods 
that provide sparse stiffmess matrices 
\RScite{shu-et-al:2005-1,sarler:2007-1,stevens-et-al:2009-1,shen:2010-1,%
yao-et-al:2012-1}.
%****************************************************************
%****************************************************************
%\section{Optimal Recovery}\RSlabel{SecOR}
%****************************************************************
%\section{Greedy Testing}\RSlabel{SecGT} 
%%%%%%%%%%%%%%%%%%%%%%%%%%%%%%%%%%%%%%
%%%%%%%%%%%%%%%%%%%%%%%%%%%%%%%%%%%%% 
\bibliographystyle{plain}
%\bibliography{RSbib,newproj}

\end{document}